\documentclass[reqno]{amsart}
\usepackage{mathtools}
\usepackage{microtype}
\usepackage{enumitem}
\usepackage{csquotes}
\usepackage{lmodern}
\usepackage[british]{babel}
\usepackage[normalem]{ulem}

\DeclareMathOperator{\curl}{curl}
\DeclareMathOperator{\Div}{div}
\DeclareMathOperator{\rmD}{dom}

\DeclarePairedDelimiter{\norm}{\lVert}{\rVert}
\DeclarePairedDelimiter{\innerproduct}{\langle}{\rangle}
\DeclarePairedDelimiter{\parentes}{\lparen}{\rparen}
\DeclarePairedDelimiter{\bracket}{\lbrack}{\rbrack}

\newcommand{\complexes}{\mathbb C}
\newcommand{\integers}{\mathbb Z}
\newcommand{\reals}{\mathbb R}
\newcommand{\naturals}{\mathbb N}

\newcommand{\qb}{\mathsf{q}_b}
\newcommand{\qmb}{\mathsf{q}^{(m)}_b}
\newcommand{\qmt}{\widetilde{\mathsf{q}}^{(m)}_b}
\newcommand{\qme}{\mathsf{q}^{(m)}_\star}
\newcommand{\Sam}{\mathcal{S}}
\newcommand{\Ham}{\mathcal{H}}
\newcommand{\Hm}{\Ham^{(m)}}

\newcommand{\Smt}{\widetilde{\Sam}^{(m)}_b}
\newcommand{\Sm}{\Sam^{(m)}_b}
\newcommand{\Sme}{\Sam^{(m)}_\star}
\newcommand{\Hb}{\Ham}
\newcommand{\Ab}{\mathbf{A}}% Magnetic vector potential
\newcommand{\M}{\mathsf{M}}% Whittaker M
\newcommand{\W}{\mathsf{W}}% Whittaker W
\newcommand{\disk}{\mathcal{B}}

\newcommand{\cG}{\mathcal{G}}
\newcommand{\cH}{\mathcal{H}}% Hamiltonian
\newcommand{\cI}{\mathcal{I}}% Momentum
\newcommand{\cL}{\mathcal{L}}% Landau level
\newcommand{\cO}{\mathcal{O}}% Ordo
\newcommand{\dd}{\mathop{}\!{\mathrm{d}}}
\newcommand{\ext}{\mathrm{ext}}
\newcommand{\Omegaext}{\Omega^{\ext}}

\newcommand{\ii}{\mathrm i}% imaginary i
\newcommand{\ee}{\mathrm e}% exponential e

\newcounter{counter_a}
\newenvironment{myenum}{\begin{list}{{\rmfamily(\roman{counter_a})}}%
{\usecounter{counter_a}
\setlength{\itemsep}{1.ex}\setlength{\topsep}{0.8ex}
\setlength{\leftmargin}{5ex}\setlength{\labelwidth}{5ex}}}{\end{list}}

\usepackage{xcolor}
\definecolor{darkgreen}{rgb}{0.0, 0.55, 0.0}

\definecolor{darkblue}{rgb}{0.0, 0.0, 0.55}
\definecolor{darkred}{rgb}{0.55, 0.0, 0.0}

\usepackage[hidelinks]{hyperref}
\hypersetup{
    colorlinks,
    citecolor=darkblue,
    filecolor=black,
    linkcolor=darkblue,
    urlcolor=black
}

\newtheorem{theorem}{Theorem}[section]

\newtheorem{lemma}[theorem]{Lemma}
\newtheorem{proposition}[theorem]{Proposition}

\newtheorem{corollary}[theorem]{Corollary}
\theoremstyle{remark}
\newtheorem{remark}[theorem]{Remark}

\numberwithin{equation}{section}

\title[Magnetic Laplacian in exterior domains]
      {On the Laplace operator\\ with a weak magnetic field\\ in exterior domains}

\author[A. Kachmar]{Ayman Kachmar}
\author[V. Lotoreichik]{Vladimir Lotoreichik}
\author[M. Sundqvist]{Mikael  Sundqvist}

\address[A. Kachmar]{The Chinese University of Hong Kong,  Shenzhen,  Guangdong,  518172,  P.R.  China.}
\email{akachmar@cuhk.edu.cn}

\address[V. Lotoreichik]{Department of Theoretical Physics, Nuclear Physics Institute, Czech Academy of Sciences, 25068, \v{R}e\v{z}, Czech Republic}
\email{lotoreichik@ujf.cas.cz}

\address[M. Sundqvist]{Department of Mathematics, Lund University, Sweden}
\email{mikael.persson\_sundqvist@math.lth.se}

\keywords{magnetic Laplacian, eigenvalue asymptotics, weak magnetic fields, isoperimetric inequality, Landau levels}
\subjclass{35P15; 58J50, 34L15, 35Q40}

\begin{document}

\begin{abstract}
  We study the magnetic Laplacian in a two-dimensional  exterior domain with
  Neumann boundary  condition and uniform magnetic field. For the exterior of
  the disk we establish accurate asymptotics of the low-lying eigenvalues in the
  weak magnetic field limit. For the exterior of a star-shaped domain, we obtain
  an asymptotic upper bound on the lowest eigenvalue in the weak field limit,
  involving the \(4\)-moment,  and  optimal for the case of the disk. Moreover,
  we prove that, for moderate magnetic fields, the exterior of the disk is a
  local maximizer for the lowest eigenvalue under a \(p\)-moment constraint.
\end{abstract}

\maketitle

\section{Introduction}

\subsection{Background and motivation} 

Inequalities between eigenvalues in terms of geometric data are central in
spectral geometry. Many interesting results are proven for the Laplace operator.
We mention the Faber--Krahn inequality~\cite{faber,krahn}, that says that among
all domains of finite volume, the ball minimizes the first Dirichlet eigenvalue.
On the contrary, the ball is a maximizer of the first non-zero eigenvalue in the
case of Neumann boundary conditions~\cite{sz,weinberger}.

L.~Erd\H{o}s extended the Faber--Krahn inequality to the magnetic Laplace
operator with a uniform magnetic field, still under a Dirichlet boundary
condition~\cite{erdos96}. A quantitative version of the isoperimetric inequality
by Erd\H{o}s  was recently obtained in~\cite{GLM23}. Much less is known for other
boundary conditions. It was suggested in~\cite{fohe19} that the disk might in
fact be a \emph{maximizer} among all bounded two-dimensional simply-connected
domains with fixed area. This conjecture is still out of reach, but there is
some  progress~\cite{clps,clps2,KL}.  Related results concern eigenvalue
optimization for the magnetic Laplacian with Robin boundary
conditions~\cite{DKL23, KL22}.

Optimization of the lowest eigenvalue for the non-magnetic Laplace operator with
attractive Robin boundary conditions in exterior domains was analysed in
~\cite{B24, KL18, KL20},   see also the monograph~\cite{BW23}.   Less is known in
this context for the magnetic Laplace operator in exterior domains, where the
optimization of the lowest eigenvalue is meaningful also for Neumann boundary
conditions.

We study here the case of magnetic Neumann boundary conditions, and weak
magnetic fields, in the exterior of star-shaped domains. The exterior of the
disk is the main suspect to be a maximizer of the lowest eigenvalue. Therefore,
we develop in Section~\ref{sec:disk} precise estimates of the eigenvalues of the
magnetic Laplacian in the exterior of a disk, and we obtain the eigenvalue
splitting in the limit of weak magnetic field. We discuss the exterior of
star-shaped domains about the origin in Section~\ref{sec:loc-opt}, and show that
the exterior of the disk is indeed a local maximizer of the lowest eigenvalue
under a fixed \(p\)-th moment of the exterior domain. Moreover, we obtain in the
same section an asymptotic upper bound on the lowest eigenvalue in the weak
field limit, which yields optimality of the exterior of a disk for sufficiently
weak magnetic fields under fixed fourth moment of the exterior domain.   

The main results in this paper are stated below in Theorem~\ref{thm:disk}
concerning the exterior of the disk, and Theorems~\ref{thm:ub-general} and
\ref{thm:loc-opt} concerning the exterior of a star-shaped domain about the
origin. They highlight two new phenomena, the first is that the weak field limit
in exterior domains turns out to have a \emph{semi-classical} character,  due to
the fact the spectrum of the non-magnetic Laplacian is purely continuous in
exterior domains.  The second is the emergence of the moment constraint under
which we verify the local optimality of the disk's exterior. A part of the
novelty in  our results, we believe that our analysis paves the way towards
understanding whether the exterior of the disk is a maximizer among all exterior
domains with fixed area or perimeter of the interior complement.

\subsection{The magnetic Laplacian in an exterior domain}

Consider a bounded simply connected domain \(\Omega\subset\reals^2\) with  a
$C^\infty$-smooth boundary \(\partial\Omega\) and the outer unit normal vector
field $\nu$. Throughout this paper,  we denote by $\Omegaext =
\reals^2\setminus\overline{\Omega}$, the exterior of the domain $\Omega$. We
also let \(\disk_R = \{x\in\reals^2\colon |x|<R\}\) denote the open disk of
radius \(R > 0\), and in particular \(\disk_1\) the unit disk. 

Let \(\Ab\) be the vector potential defined as
\begin{equation}\label{eq:def-A}
  \Ab(x) = \frac12(-x_2,x_1)^\top,
\end{equation}
and which generates the  unit magnetic field, \(\curl\Ab=1\). We study the
self-adjoint magnetic Neumann Laplacian in $L^2(\Omegaext)$,
\[
\begin{split}
  \Hb u &:= (-\ii\nabla -b\Ab)^2 u,\\
  	\rmD(\cH) &:= 
  	\big\{u\in H^1_{\Ab}(\Omegaext)\colon (-\ii\nabla - b\Ab)^2u\in L^2(\Omegaext),\, 
  	\nu\cdot(\nabla-\ii b\Ab) u = 0~\text{on}~\partial\Omega \big\},
\end{split}
\]
where  \(b \in \reals\) and the first-order magnetic Sobolev space is defined by
\(H^1_{\Ab}(\Omegaext) := \{u\in L^2(\Omegaext)\colon (-\ii\nabla-b\Ab)u \in
 [L^2(\Omegaext)]^2\}\). In fact, the operator is unitarily equivalent 
under a sign-flip of the magnetic field, so it is sufficient to consider \(b \geq 0\). The
operator \(\Hb\) is associated with the quadratic form
\[
  \qb(u):=\int_{\Omegaext} |(-\ii\nabla - b\Ab)u|^2\dd x,
  \qquad 
  \rmD(\qb) := H^1_{\Ab}(\Omegaext).
\]
Contrary to the case of an interior domain, the spectrum of $\Hb$ depends on the
gauge of the vector potential $\Ab$, since we can not rule out this dependence
by a gauge transformation.

The essential spectrum of \(\Hb\) consists of the Landau levels,
\(\{(2n-1)b\colon n\in\naturals\}\), and the discrete spectrum is infinite and
accumulates at each Landau level from below~\cite{GKS}. In particular, the
spectrum below $b$ consists of a sequence of discrete eigenvalues
\[  
  \lambda_1(b,\Omegaext)
  \leq
  \lambda_2(b,\Omegaext)
  \leq 
  \ldots 
  <
  b,
\]
counted with multiplicity. 

The eigenvalues can, as usual, be characterized by the variational min-max theorem,
\[ 
  \lambda_k(b,\Omegaext)
  =
  \inf_{ \substack{M\subset H^1_{\Ab}(\Omegaext) \\ \dim M = k} }
  \biggl(
    \max_{\substack{u\in M\\u\neq 0}} \frac{\qb(u)}{\|u\|_{L^2(\Omegaext)}^2}
  \biggr).
\]

\subsection{Exterior of the disk}

We will concentrate on the exterior of the disk, and prove that for weak
magnetic fields, all the eigenvalues below $b$ are simple, and obey a two-term
asymptotics.

\begin{theorem}\label{thm:disk} 
  Assume that \(\Omega = \disk_R\) and that \(b > 0\). Then the ground state of
  \(\Hb\) will be non-radial. Moreover, for any $k\geq 1$,
  \[
    \lambda_k(b,\disk_R^{\ext}) 
    = 
    b - \frac{R^{2k}}{2^{k - 1}(k - 1)!}b^{k + 1} + \cO(b^{k + \frac32}),
    \quad (b \to 0^+).
  \]
\end{theorem}

\begin{remark}\label{rem:thm-disk}
\begin{myenum}
  \item By a scaling argument, it suffices to prove Theorem~\ref{thm:disk} in
  the case of the exterior of the unit disk, \(\Omega = \disk_1\), which
  corresponds to \(R = 1\).
  \item  Our proof of Theorem~\ref{thm:disk} relies on separation of variables
  and the analysis of fibre operators. Even though this method does not directly
  carry over to general domains, we will give in Theorem~\ref{thm:ub-general} an
  asymptotic upper bound on the lowest eigenvalue for star-shaped domains that
  will enable us to compare with the disk under a certain moment constraint.
  \item For the Neumann magnetic operator in the disk \( \disk_R
  \), ground states are radial at least for \(b R^2<1\)
  (see~\cite[Prop.~2.3]{KL}). The limit of the lowest eigenvalue as \(b \to
  0^+\) follows from standard perturbation theory~\cite[Prop. 1.5.2]{FH-b}; in
  fact, as \(b \to 0^+\),
  \[
    \lambda_1(b,\disk_R) =  \frac{b^2R^2}{8} + \cO(b^3).
  \]
  \item Interestingly,  there is an  analogy between Theorem~\ref{thm:disk} and
  the strong field asymptotics for the Dirichlet eigenvalues in the interior of
  the disk,  which highlights an interior-exterior,   Dirichlet-Neumann and
  weak-strong field duality.   In fact,~\cite[Thm.~5.2]{HS} yields that,
  \[
    \lambda_k^{\mathrm{D}}(b,\disk_R)
    =
    b
    +
    \frac{R^{2k}}{2^{k - 1}(k - 1)!}b^{k + 1}\ee^{-bR^2/2}
    +
    o\bigl( b^{k + 1}\ee^{-bR^2/2}\bigr),
    \quad
    (b \to +\infty),
  \]
  where $\lambda_k^{\mathrm D}(b,\disk_R)$ denotes the $k$-th eigenvalue of
  $(-\ii\nabla-b\Ab)^2$ in $L^2(\disk_R)$,  with Dirichlet boundary condition on
  $\partial\disk_R$.
\end{myenum}
\end{remark}

\subsection{Exterior of a star-shaped domain.}

In our considerations of general domains  we employ a quantity associated
merely to the exterior domain itself. Clearly, the area of the exterior domain
is always infinite, and we cannot use it in the analysis. However, the
following $p$-th moment
\[
  \cI_p(\Omegaext ) := \int_{\Omegaext }\frac{1}{|x|^p}\dd x,\qquad p > 2,
\]
is finite, provided that the origin belongs to $\Omega$. Such a moment appears
in the geometric constraint and in the asymptotic upper bound for the values of
$p$ specified below. In some cases, suitable values of $p$ will also depend on
the intensity of the magnetic field $b$. We also remark that the $p$-th moment
of $\Omegaext$ is not invariant under translations of the domain and the results
below can be further optimized by choosing the origin so that $\cI_p(\Omegaext)$
is minimal. 

\subsubsection{Asymptotic upper bound.}
In the case where  the domain $\Omega$
is star-shaped with respect to the origin and parametrized by a smooth mapping
$\reals/[0,2\pi]\ni\theta\mapsto\rho(\theta)$ as $\Omega = \{(r,\theta)\colon r
<\rho(\theta)\}$,  we establish an upper bound on the lowest eigenvalue $\lambda_1(b,\Omegaext)$,  valid in the limit of small $b$.
\begin{theorem}\label{thm:ub-general}
Let $\Omega \subset\reals^2$ be a  bounded $C^\infty$-smooth domain
star-shaped with respect to the origin and parametrized by the $C^\infty$-smooth
function $\rho\colon \reals/[0,2\pi]\to \reals_+$. Then, as $b\to 0^+$,  we have
\[ 
  \lambda_1(b,\Omegaext)
  \leq 
  b
  -
  \Bigl(
    \frac{|\Omega|}{2\pi}+\frac{\pi}{2\cI_4(\Omegaext)} 
  \Bigr)
  b^2 
  + 
  \cO(b^{5/2}).
\]
Moreover, for $\Omega$ being not a disk,  there exists $b_0=b_0(\Omega)>0$ such that,  for $0<b<b_0$,  we have, 
\[
  \lambda_1(b,\Omegaext)
  <
  \lambda_1(b,\disk_{R_\star}^{\ext}),
\]
where $R_\star$ denotes the radius of the disk satisfying
$\cI_4(\disk_{R_\star}^{\ext})=\cI_4(\Omegaext)$. 
\end{theorem}

\begin{remark}\label{rem:thm-ub-general}
  \begin{myenum}
  \item  The upper bound in Theorem~\ref{thm:ub-general} is consistent with the
  asymptotics in Theorem~\ref{thm:disk} if we take $\Omega=\disk_R$,   the disk
  of radius $R$.
  \item If we choose the radius $R$ so that $|\disk_{R}|=|\Omega|$,   the upper
  bound in Theorem~\ref{thm:ub-general} does not yield that the exterior of the
  disk is a maximizer of the lowest eigenvalue for sufficiently small $b$.
  \end{myenum}
\end{remark}
\subsubsection{Moderate field and moment constraint.}
We consider  next the more general  case with a moderate magnetic field and
address whether the exterior of the disk is a local maximizer of the lowest
magnetic Neumann eigenvalue under suitable geometric constraints.   We provide
an affirmative answer when imposing a certain  constraint and when $\Omega$ is a
small deformation of the disk in a certain sense,  which will be clarified
below. 

Our result demonstrates that the exterior of the disk is a local optimizer under
the $p$-moment constraint. In the formulation of the theorem below, we use    
an auxiliary quantity
\[
  p_\star(b) \coloneqq 2 + \frac{16b - 4b^2}{4 - 8b + b^2}.
\]
\begin{theorem}\label{thm:loc-opt}
	Let $R>0$, $bR^2 \in (0,4-2\sqrt{3})$ and $p > p_\star(bR^2)$. Let $\Omega
	\subset\reals^2$ be a bounded $C^\infty$-smooth domain star-shaped
	with respect to the origin and parametrized by the $C^\infty$-smooth function
	$\rho\colon \reals/[0,2\pi]\to \reals_+$. Then, there exists $r_\star =
	r_\star(b,p,R) \in (0,R)$ such that under the assumptions $\cI_p(\Omegaext ) =
	\cI_p(\disk_{R}^{\ext})$ and $R - r_\star \le \rho \le R + r_\star$, the
	following isoperimetric inequality holds:
	\[
		\lambda_1(b,\Omegaext ) \le \lambda_1(b,\disk_{R}^{\ext}).
	\]
\end{theorem}

\begin{remark}\label{rem:opt}
  Natural constraints under which optimality of the exterior of the disk can be
  analysed are fixed area or fixed perimeter of the interior domain, as in
  related analysis on the optimization of the lowest eigenvalue for the Robin
  Laplacian on exterior domains~\cite{B24, KL18, KL20}. The method we propose
  for the magnetic Neumann Laplacian does not imply local optimality of the
  exterior disk under such constraints. (Local) optimality of the exterior disk
  under these constraints remains an open problem. The magnetic problem is
  subtle due to the fact that the ground state is not a radial function and the
  problem shares some common features with the optimization of the second Robin
  eigenvalue in exterior domains~\cite{KL23}. 
\end{remark}

\begin{remark}
	\begin{myenum}
	\item In the limit $b\to 0^+$ we have $p_\star(b)\to 2$. Thus, for any $p > 2$,
	we get by Theorem~\ref{thm:loc-opt} local optimality of the exterior of the disk for a fixed $p$-moment provided that $b > 0$ is sufficiently small. 
	\item We also point out that Theorem~\ref{thm:loc-opt} establishes uniform local optimality of the exterior of
	the disk
	(under fixed $p$-moment)
	within exteriors of all smooth star-shaped domains, whose Hausdorff distance from the disk is smaller than some fixed positive constant, which depends only on the intensity $b$ of the magnetic field, the parameter $p$, and the radius $R$ of the disk. 
\end{myenum}	
\end{remark}

We organized the rest of the paper as follows. Section~\ref{sec:disk} is devoted
to the proof of Theorem~\ref{thm:disk}. Section~\ref{sec:loc-opt} is devoted to
the proof of Theorems~\ref{thm:ub-general} and \ref{thm:loc-opt}. Finally, Appendix~\ref{sec:app} contains
a result on the intersection of eigenvalue curves of the fibre operators.

\section{The weak field limit in the exterior of the disk}\label{sec:disk}

In this section we prove the two-term asymptotic formula of the smallest
eigenvalue in the exterior of the disk, stated in Theorem~\ref{thm:disk}.

\subsection{Fibre operators in the exterior of the disk}

\subsubsection{Splitting into fibre operators}
We work in the exterior of the unit disk, \(\disk_1^{\ext}\). The action of the
operator $\cH$ reads in polar coordinates as,
\[
  \Hb = -\frac{\partial^2}{\partial r^2} 
        - \frac{1}{r}\frac{\partial}{\partial r} 
        + \Bigl(\frac{-\ii\partial_\theta}{r} - \frac{b r}{2}\Bigr)^2.
\]
Since \(\{\ee^{\ii m\theta}\}_{m\in\mathbb{Z}}\) is an orthogonal basis for the
square integrable functions on the circle, our operator splits naturally. We can
work with the family \(\{\Hm\}_{m\in\integers}\) of self-adjoint fibre operators
in \(L^2((1,+\infty);r\dd r)\),
\begin{equation}\label{eq:Hm-b}
\begin{split}
	\Hm  &= \Hm(b) := -\frac{\dd^2}{\dd r^2} -\frac1r\frac{\dd}{\dd r} + V_{m,b}(r),\\
	\rmD(\Hm) &:= \left\{u\colon u,-u''- \frac{u'}{r} + V_{m,b}u \in L^2((1,+\infty);r\dd r),\,\, u'(1) = 0\right\},
\end{split}	
\end{equation}
where
\[
  V_{m,b}(r) = \left(\frac{m}{r}-\frac{br}{2}\right)^2.
\]
For each
\(b\in\mathbb{R}\),
\[
  \sigma(\Hb(b)) = \overline{\bigcup_{m\in\mathbb{Z}} \sigma(\Hm(b))},
\]
where we use the decomposition $\cH\simeq\bigoplus_{m\in\mathbb{Z}}\Hm$. We let
\(\lambda_{1}^{(m)}(b) \leq \lambda_{2}^{(m)}(b) \leq \ldots\) denote the
eigenvalues of \(\Hm\). Away from eventual crossings, all \(b \mapsto
\lambda_k^{(m)}(b)\) are real analytic,  see \cite[Chap.~VII, \S 3.5, Thm.\
3.9]{K}.

We remark that the ground state energy of $\Hb$ can be expressed as
\begin{equation}
  \lambda_1(b,\disk_1^{\ext}) = \inf_{m \in \integers}\lambda_1^{(m)}(b).
\end{equation}
In fact, we know that the essential spectrum of \(\Hb\) consists of the Landau
levels \(b\), \(3b\), and so on. We will find eigenvalues of \(\Hm(b)\) that are
less than \(b\), so \(\lambda_1(b,\disk_1^{\ext})\) will indeed be an eigenvalue.

\subsubsection{A first control of each eigenvalue}

We recall that we only consider positive \(b\), and warm up by excluding all \(m
\leq 0\) from the game. We also show that each \(b \mapsto \lambda_1^{(m)}(b)\)
will cross the lowest Landau level \(b\) exactly once for \(b > 0\), and that it
will happen for \(b = 2m\).

\begin{proposition}\label{prop:disk-m}
For \(b > 0\) the following holds.
\begin{enumerate}[wide,leftmargin=0pt,label=\emph{(\arabic*)}]
  \item If \(m < 0\), then \(\lambda_1^{(m)}(b) \geq 2b\).
  \item If \(m = 0\) and \(b > 0\), then \(\lambda_1^{(m)}(b) > b\).
  \item If \(m\geq 1\) then
    \[ 
        \begin{cases}
          \lambda_1^{(m)}(b) < b & \text{if }b < 2m,\\
          \lambda_1^{(m)}(b) = b & \text{if }b = 2m,\\
          \lambda_1^{(m)}(b) > b & \text{if }b > 2m.
        \end{cases}
    \]
  \item If \(m\geq 1\), then the second eigenvalue of \(\cH^{(m)}\) satisfies
  \(\lambda_2^{(m)}(b) > b\).
\end{enumerate}
\end{proposition}

\begin{proof}
We show the statements one by one.
\begin{enumerate}[wide,leftmargin=0pt]
  \item For \(m < 0\), we  have
  \[
    V_{m,b}(r) \geq 2|m|b
  \]
  which directly yields
  \[
    \lambda_1^{(m)}(b) > 2 b.
  \]
  \item For \(m = 0\) and \(b > 0\), the eigenvalue equation \(\Ham^{(0)} u =
  b u\) has the general solution
  \[
    u(r)
    =
    u_0(r) 
    \coloneqq 
    \ee^{-br^2/4} 
    \left(
      c_1 + c_2 \int_1^r \rho^{-1} \ee^{b\rho^2/2}\dd\rho
    \right).
  \]
  For the function not to increase exponentially at \(+\infty\) we must have
  \(c_2 = 0\). But then we find that \(u_0'(r) = -c_1(br/2)\ee^{-br^2/4}\), and
  so we cannot satisfy the Neumann boundary condition \(u_0'(1) = 0\) unless
  also \(c_1 = 0\). Hence, there is no \(b > 0\) such that \(\lambda_1^{(0)}(b)
  = b\). We notice next that \(V_{0,b}(r) = b^2r^2/4 \geq b^2/4\), so if \(b >
  2\) it holds that \(\lambda_1^{(0)}(b) > b\). By continuity of \(b \mapsto
  \lambda_1^{(0)}(b)\), this inequality is valid for all \(b > 0\).

  \item Let us fix \(m \geq 1\), solving again the eigenvalue equation for the
  lowest Landau level, \(\Hm u = b u\), we get
  \[
    u (r) 
    = 
    u_m (r)
    \coloneqq r^m \ee^{-br^2/4} 
    \left(
      c_1 + c_2\int_1^{r}\rho^{-1-2m}\ee^{b\rho^2/2}\dd\rho 
    \right).
  \] 
  The term with the coefficient \(c_2\) is not in \(L^2((1,+\infty);r\dd r)\), so
  we take \(c_2 = 0\). We also choose \(c_1\) so that \(u_m\) becomes
  normalized. The solution \(u_m\) satisfies the Neumann boundary condition if
  and  only if \(b = 2m\). Thus,
  \begin{equation}\label{eq:b=2m}
    b \in \sigma(\Hm) \iff b = 2m. 
  \end{equation}
  Moreover, since \(V_{m,b}(r) \geq b^2/4 - mb\) it holds that \(V_{m,b}(r) >
  b\) if \(b\) is sufficiently large. This implies that \(\lambda_1^{(m)}(b) >
  b\) if \(b\) is sufficiently large. We conclude that \(\lambda_1^{(m)}(b) =
  2m\) precisely when \(b = 2m\) and that \(\lambda_1^{(m)}(b) > b\) for all \(b
  > 2m\). 

  We next show that \((\lambda_1^{(m)})'(2m) > 1\); that implies that we really
  have a crossing at \(b = 2m\), and that \(\lambda_1^{(m)}(b) < b\) for \(0 < b
  < 2m\). By the Feynmann--Hellman formula
   \[
    (\lambda_1^{(m)})'(b) = -\int_1^{+\infty} \left(m-\frac{br^2}{2}\right)|u_m|^2r\dd r,
  \]
  which can be justified by the fact that the sesquilinear forms associated with
  $\Hm(b)$ form a holomorphic family of the type (B) in the sense of Kato with
  respect to the parameter $b$ and thus the perturbation formulae
  in~\cite[Chap.~VII, \S 4.6]{K} apply. Since \(1 = \int_1^{+\infty}|u_m|^2 r\dd
  r\), we get for \(b = 2m\),
  \[
    (\lambda_1^{(m)})'(2m) - 1 
    = 
    c_1^2 \int_1^{+\infty} \big(m(r^2-1)-1\big) r^{2m+1} \ee^{-mr^2} \dd r.
  \]
  Integrating by parts, we get
  \[
    \begin{aligned}
      \int_1^{+\infty}  m(r^2-1)r^{2m+1}\ee^{-mr^2} \dd r
      & = \int_1^{+\infty}
        \left(
          1 + m - \frac{m}{r^2}
        \right)
        r^{2m+1}\ee^{-mr^2}\dd r\\
      & > \int_1^{+\infty}r^{2m+1}\ee^{-mr^2}\dd r,
    \end{aligned}
  \]
  which yields that \((\lambda_1^{(m)})'(2m) > 1\).
  \item
  By the previous step we already know that \(\lambda_2^{(m)}(b) > b\) for \(b >
  2m\). The inequality must still hold for \(b = 2m\) since the groundstate of
  \(\Hm\) is simple. In fact, it must also hold for \(0 < b < 2m\), since
  otherwise we would have \(\lambda_2^{(m)}(b) = b\) for some \(0 < b < 2m\),
  and that is not possible according to~\eqref{eq:b=2m}.\qedhere
\end{enumerate}
\end{proof}

\begin{corollary}\label{prop:disk} 
  For all \(b > 0\) we have \(\lambda_1^{(m)}(b) < b\) and \(\Hb\) does not have
  radially symmetric ground states.
\end{corollary}

\begin{proof}
  By Proposition~\ref{prop:disk-m}, the ground state must be of the form
  \(u_m(r)\ee^{\ii m\theta}\) with \(m\) a positive integer.
\end{proof}

\begin{corollary}\label{prop:disk*} 
  For all \(0<b<2\) and $m\geq 1$, $\lambda_1^{(m)}(b)<b$ is a discrete
  eigenvalue of $\Hb$; moreover, $\lambda_1^{(m)}(b)\to b$ as $m\to +\infty$.
\end{corollary}
\begin{proof}
  Proposition~\ref{prop:disk-m} tells us that
  \(\{\lambda_1^{(m)}(b)\}_{m\geq1}\) is the sequence of discrete eigenvalues of
  \(\Hb\) below the bottom of the essential spectrum, \(b\). Hence,
  \(\lambda_1^{(m)}(b)\to b\) as \(m\to +\infty\). 
\end{proof}

\subsection{An effective operator}

We would like to understand the limit behaviour,  as $b \to 0^+$,  of the
eigenvalues $\lambda_1^{(m)}(b)$ and $\lambda_2^{(m)}(b)$ of the operator
$\Hm(b)$ introduced in~\eqref{eq:Hm-b}. This will be done by deriving an
effective operator, which allows us to prove the following proposition.

\begin{proposition}\label{prop:spectrum-1-2}
For every $m\geq 1$, we have as $b\to 0^+$,
\[
  \lambda_1^{(m)}(b)=b+o(b),\quad \lambda_2^{(m)}(b)\geq 3b+o(b).
\]
\end{proposition} 

The rest of this subsection is devoted to the proof of
Proposition~\ref{prop:spectrum-1-2}.

\subsubsection{Unitary transformation and scaling} 

In order to transform the weighted space $L^2((1,+\infty);r\dd  r)$ to the
non-weighted Hilbert space $L^2((1,+\infty);\dd  r)$, we apply the unitary
transformation $U\colon L^2((1,+\infty);\dd r) \to L^2((1,+\infty); r\dd r)$
\[
  U\colon u(r) \mapsto r^{-1/2}u(r),
\]
and then, we do the change of variable
\[
  r\mapsto \sqrt{\frac{b}{2}}\,(r-1)
\]
and get
\[
  U^{-1}\Hm(b)U\simeq\frac{b}{2}\,\Sm,
\]
where the self-adjoint operator $\Sm$ in $L^2(\mathbb{R}_+)$ is given by
\[
  \begin{aligned}
    \Sm
    & =
    -\frac{\dd^2}{\dd r^2}+w_{m,b}(r),\\
    \rmD(\Sm) &= \left\{u\colon u,-u'' +w_{m,b} u\in L^2(\mathbb{R}_+),\,\, u'(0)
    =
    \sqrt{\frac{1}{2b}}\,u(0)\right\},
  \end{aligned}
\]
with the potential $w_{m,b}$  defined by
\[
  w_{m,b}
  \coloneqq
  \frac{4m^2-1}{4\left(r + \sqrt{\frac{b}{2}}\right)^2}
  + \left(
      r+\sqrt{\frac{b}{2}}
    \right)^2-2m.
\]
\subsubsection{Localizing the spectrum of $\Sm$}
For fixed $m$, the spectrum of $\Sm$ converges to that of an operator on
$\mathbb{R}_+$ formally defined by
\[
  \Sme
  =
  -\frac{\dd^2}{\dd r^2}+\frac{4m^2-1}{4r^2}+r^2-2m,
\]
with the boundary condition $u(0)=0$; (we introduce $\Sme$ rigorously as the
self-adjoint operator associated with the closed, densely defined, non-negative
quadratic form $\qme$ defined below). We will make this convergence precise in
Proposition~\ref{prop:eff*}, but let us first  introduce the two eigenvalues
\begin{equation}%\label{eq:def-ev}
  \begin{aligned}
    \mu_1(\Sm)
    &=
    \inf_{u\in D(\qmb)\setminus\{0\}}\frac{\qmb(u)}{\|u\|^2_{L^2(\reals_+)}},\\
    \mu_1(\Sme)
    &=
    \inf_{u\in D(\qme)\setminus\{0\}}\frac{\qme(u)}{\|u\|^2_{L^2(\reals_+)}},
  \end{aligned}
\end{equation}
where
\[
  \begin{aligned}
    \qmb(u)
    &=
    \int_{\reals_+}\left(|u'(r)|^2+ w_{m,b}(r)|u(r)|^2\right)\dd r
    +\sqrt{\frac1{2b}}\,|u(0)|^2,\\
    \qme(u)
    &=
    \int_{\reals_+}\left(|u'(r)|^2+w_{m,0}(r)|u(r)|^2 \right)\dd r,
  \end{aligned}
\]
and\footnote{By the Hardy inequality,  $\int_{\reals_+} r^{-2}|u|^2\dd r\leq 4\int_{\reals_+}|u'|^2\dd r$ for all $u\in H^1_0(\reals_+)$, so we get for free that $u/r\in L^2(\reals_+)$.}
\[
  \begin{aligned}
    \rmD(\qmb) &= \{u\in H^1(\reals_+)\colon ru\in L^2(\reals_+)\},\\
    \rmD(\qme) &= \{u\in H^1_0(\reals_+)\colon ru\in L^2(\reals_+)\}.
  \end{aligned}
\]
We also denote by $\{\mu_k(\Sm)\}_{k\geq 1}$ and $\{\mu_k(\Sme)\}_{k\geq 1}$  the sequences of eigenvalues of $\Sm$ and $\Sme$ respectively,   numbered in increasing order including multiplicity.    Observe that, the scaling and unitary transformation we did,  yield the following
\[
  \lambda_k^{(m)}(b)
  =
  \frac{b}2\mu_k(\Sm)\qquad (k\geq 1).
\]
Therefore, we have as immediate consequence of Proposition~\ref{prop:disk-m},
\begin{equation}\label{eq:ub}
  \forall\,m\geq1, \forall\,b\in(0,2),\quad \mu_1(\Sm)<2.
\end{equation}

\begin{proposition}\label{prop:eff*}
  For all $m\ge 1$, the operator $\Sm$ converges in the strong resolvent sense
  to $\Sme$ as $b\to 0^+$.  Moreover, the eigenvalues of $\Sm$ converge to the
  corresponding eigenvalues of $\Sme$ as $b\to 0^+$.
\end{proposition}

\begin{proof}
   To be able to apply the abstract result~\cite[Theorem S.14]{RS-I} we
  need a monotone family of quadratic forms. 
	Let us therefore introduce an intermediate self-adjoint  operator
	$\Smt$ in $L^2(\reals_+)$ associated with the closed, semi-bounded, symmetric and densely defined quadratic form
	\[
    \begin{aligned}
      \qmt(u) 
      &\coloneqq 
      \int_{\reals_+}
      \left(
        |u'(r)|^2 +w_{m,b}^{\mathrm{int}}(r)|u(r)|^2
      \right)
      \dd r +\sqrt{\frac{1}{2b}}|u(0)|^2,\\
      \rmD(\qmt)
      &\coloneqq 
      \{u\in H^1(\reals_+)\colon r u\in L^2(\reals_+)\},
    \end{aligned}	
	\]
	where the potential is defined by
	\[
  	\widetilde{w}_{m,b}(r)
    \coloneqq 
    \frac{4m^2-1}{4\left(r+\sqrt{\frac{b}{2}}\right)^2} + r^2 -2m.
	\]
	It is straightforward to see that the family of quadratic forms $\reals_+\ni
	b\mapsto \qmt$ is monotonously decreasing in $b$ in the sense of ordering of
	the quadratic forms. It follows from~\cite[Theorem S.14]{RS-I} using the
	identities
	\[
		 D(\qme)
		 = \left\{u\in H^1(\reals_+)\colon 	r u \in L^2(\reals_+), ~
		\sup_{b > 0} \qmt(u) < \infty\right\}
	\]
	and
	\[
		\lim_{b\to 0^+} \qmt(u) = \qme(u),
    \qquad \text{for all } u\in D(\qme),
	\]
	that the operators $\Smt$ converge to $\Sme$ in the strong resolvent sense as
	$b\to 0^+$. Moreover, for any  $b \in (0,8)$ and all $u\in D(\qmb) = D(\qmt)$,
	we have
	\[
	\begin{aligned}
		\left|\qmb(u) - \qmt(u)\right|
		&=\int_{\reals_+}\left(\sqrt{2b} r + \frac{b}{2}\right)|u(r)|^2\dd r\\
		&\le 
		\sqrt{2b}\int_{\reals_+}\left( r + 1\right)|u(r)|^2\dd r\\
    &\le
		2\sqrt{2b}\left(\qmt(u)+ (2m+1)\|u\|^2_{L^2(\reals_+)}\right).\\
	\end{aligned}
	\]
	Hence, we get by~\cite[Theorem VI 3.4]{K} that
	\begin{equation}\label{eq:normresolvent}
		\big\|(\Sm + 2m+ 1)^{-1} - (\Smt+ 2m+1)^{-1}\big\|
    \to 0
    \quad\text{ as } b\to 0^+.
	\end{equation}
	It follows from the fact that $\Smt$ converge to $\Sme$ in the strong
	resolvent sense together with~\eqref{eq:normresolvent} that $\Sm$ also
	converge to $\Sme$ in the strong resolvent sense as $b\to 0^+$.
	
	Note that $\Sm \ge \mathcal T^{(m)}$ where the semi-bounded self-adjoint
	operator $\mathcal T^{(m)}$ in $L^2(\reals_+)$ with purely discrete spectrum,
	associated with the form
	\[
		\{u\in H^1(\reals_+)\colon r u\in L^2(\reals_+)\}\ni u\mapsto
		\int_{\reals_+}(|u'(r)|^2 + (r^2-2m)|u(r)|^2)\dd r.
	\]
	Thus, we get by~\cite{W80} that the eigenvalues of $\Sm$ converge to the
	corresponding eigenvalues of~$\Sme$.	 
\end{proof}	

\subsubsection{Spectral gap for the effective operator}

Now we turn to calculating the lowest eigenvalue of $\Sme$ and to estimating the
spectral gap between the first and second eigenvalues.   We will compare with the
Landau Hamiltonian to obtain the gap.

\begin{proposition}\label{prop:eff0}
  For all $m\geq 1$, we have $\mu_1(\Sme)= 2$ and $\mu_2(\Sme)\geq 6$.
\end{proposition}

Note that Proposition~\ref{prop:eff0} yields the conclusion in
Proposition~\ref{prop:spectrum-1-2}, thanks to Proposition~\ref{prop:eff*}. The
proof of Proposition~\ref{prop:eff0} relies on the following observation.

\begin{lemma}\label{lem:identity} 
  Fix a non-negative integer $m$ and suppose that $g\in \rmD(\Sme)$,
  with $g\neq 0$. Let $h\in L^2(\reals^2)$ be defined in polar coordinates as
  \[
    h=r^{-1/2}g(r)\ee^{\ii m\theta}.
  \]
  Then, 
  \[
    h\in H^1_{2\Ab}(\reals^2)
    \coloneqq
    \{u\in L^2(\reals^2)\colon (-\ii \nabla-2\Ab)u\in L^2(\reals^2)\},
  \]
  and
  \begin{equation}\label{eq:eff-lb-1}
    \frac{\innerproduct{\Sme g,g}_{L^2(\reals_+)}}
         {\|g\|^2_{L^2(\reals_+)}}
    =
    \frac{\|(-\ii\nabla-2\Ab)h\|^2_{L^2(\reals^2)}}
         {\|h\|^2_{L^2(\reals^2)}}.
  \end{equation}
  Moreover, if $g$ is orthogonal to $u$ in $L^2(\reals_+)$, with
  $u(r)=r^{m+1/2}\ee^{-r^2/2}$, then
  \begin{equation}\label{eq:eff-lb-2}
    \frac{\|(-\ii\nabla-2\Ab)h\|^2_{L^2(\reals^2)}}
         {\|h\|^2_{L^2(\reals^2)}}
    \geq 6.
  \end{equation}
\end{lemma}

\begin{proof}
Let $\dot\Omega=\reals^2\setminus\{0\}$ and consider the magnetic Sobolev space
\[
  H^1_{2\Ab}(\dot\Omega)
  = 
  \{u\in L^2(\dot\Omega)\colon (-\ii \nabla-2\Ab)u\in L^2(\dot\Omega)\},
\]
where the derivative is seen as a distribution on $\dot\Omega$, i.e., $(-\ii
\nabla-2\Ab)u\in \mathcal D'(\dot\Omega)$. Since $\Ab\in L^2_{\mathrm{loc}}
(\reals^2)$, we have that $H^1_{2\Ab}(\dot\Omega) = H^1_{2\Ab}(\reals^2)$
(see~\cite[Rem.\ 2.2]{KP}).

By a straightforward calculation, we have 
\[
  r^{-1/2}\Sme r^{1/2}
  =
  - \frac{\dd^2}{\dd r^2}
  - \frac{1}{r}\frac{\dd}{\dd r}
  + \left(\frac{m}{r} - r \right)^2
\]
and
\[
  \left(
    - \frac{\dd^2}{\dd r^2}
    - \frac{1}{r}\frac{\dd}{\dd r}
    + \left(\frac{m}{r} - r \right)^2
  \right)
  \psi(r)
  =
  \ee^{-\ii m\theta}
  (-\ii\nabla-2\Ab)^2
  \ee^{\ii m\theta}
  \psi(r).
\]
Writing \(g=r^{1/2}\psi\), we get  $h=\ee^{\ii m\theta}\psi$,
$(-\ii\nabla-2\Ab)^2h\in L^2(\dot\Omega)$, and
\[ 
  \begin{aligned}
    \innerproduct{\Sme g,g}_{L^2(\reals_+;\dd r)} 
    &=  
    \innerproduct{ r^{-1/2}\Sme r^{1/2}\psi,\psi }_{L^2(\reals_+;r\dd r)}\\
    &=
    \innerproduct{ (-\ii\nabla-2\Ab)^2h,h }_{L^2(\reals_+;r\dd r)}.
  \end{aligned}
\]
Integrating with respect to \(\theta\) on \([0,2\pi)\), we get
\[
  2\pi\innerproduct{ \Sme g,g }_{L^2(\reals_+;\dd r)}
  =
  \innerproduct{ (-\ii\nabla-2\Ab)^2h,h }_{L^2(\reals^2;\dd x)}
  =
  \norm{(-\ii\nabla-2\Ab)h}^2_{L^2(\dot\Omega)}.
\]
This proves that $h\in H^1_{2\Ab}(\dot\Omega)=H^1_{2\Ab}(\reals^2)$, and
yields the identity in~\eqref{eq:eff-lb-1}.

Now, we assume that the orthogonality condition on \(g\) holds, 
\begin{equation}\label{eq:eff-orth}
  \int_{\reals_+} g(r) r^{m+\frac12}\ee^{-r^2/2}\dd r = 0.
\end{equation}
Recall a characterization of the lowest Landau level
\[
  \begin{aligned}
    \cL
    &:=
    \{u\in L^2(\reals^2)\colon (-\ii\nabla -2\Ab)^2u=2 u\}\\
    &=
    \{f(z)\ee^{-|z|^2/2}\colon f \text{ is holomorphic in }\complexes\},
  \end{aligned}
\]
and that 
\[
  \inf_{\substack{v\in \cL^\perp\cap H^1_{2\Ab}(\reals^2)\\v\neq 0}}
  \frac{\|(-\ii\nabla-2\Ab)v\|^2_{L^2(\reals^2)}}
       {\|v\|_{L^2(\reals^2)}}
  = 6.
\]
Therefore, to prove~\eqref{eq:eff-lb-2}, it suffices to verify that $h\in
\cL^\perp$. Let $f$ be a holomorphic function, and notice that, for $r>0$, we
have by Cauchy's integral formula,
\[ 
  \int_0^{2\pi} r^{-m}\ee^{-\ii m\theta} f(r\ee^{\ii\theta}) \dd\theta
  =
  \frac{1}{\ii} \int_{C_r} \frac{f(z)}{z^{m+1}}\dd z
  =
  \frac{2\pi}{m!}f^{(m)}(0),
\] 
where $C_r$ is the circle centred at the origin and of radius $r$.
Consequently, with $u=f(z)\ee^{-|z|^2/2}$, we have
\[
  \begin{aligned}
    \innerproduct{ u, h }_{L^2(\reals^2)}
    &
    =
    \int_0^{+\infty}
      \left(
        \int_0^{2\pi} r^{-m}\ee^{-\ii m\theta}  f(r\ee^{\ii\theta}) \dd\theta
      \right) 
      r^{m+\frac12}\ee^{-r^2/2}\overline{g(r)}\dd r\\
    & 
    =
    \frac{2\pi}{m!}f^{(m)}(0)
    \int_{\reals_+} r^{m+\frac12}\ee^{-r^2/2}\overline{g(r)} \dd r
    =
    0,
  \end{aligned}
\]
where we used~\eqref{eq:eff-orth} in the last step.
\end{proof}

\begin{proof}[Proof of Proposition~\ref{prop:eff0}]
It is easy to check that \(u=r^{m+\frac12}\exp\left(-r^2/2\right)\) is an
eigenfunction of \(\Sme\) with eigenvalue $2$, so \(\mu_1(\Sme)\leq 2\).

By Lemma~\ref{lem:identity} and the min-max principle,  we have
\(\mu_1(\Sme)\geq 2\), and the function $u$ is a ground state of $\Sme$.  If
$g\in \mathrm{Dom}(\Sme)$ is orthogonal to $u$,  we have by~\eqref{eq:eff-lb-2},
\(\innerproduct{ \Sme g,g}_{L^2(\reals_+)}\geq 6\|g\|^2_{L^2(\reals_+)}\), and
by the min-max principle,  we get that \(\mu_2(\Sme)\geq 6\).
\end{proof}

\subsection{Refined eigenvalue asymptotics}

We improve the leading term asymptotics in Proposition~\ref{prop:spectrum-1-2}
by calculating the subleading term. 

\begin{theorem}\label{thm:lambdamasymptotic}
Let \(m \geq 1\). As \(b\to 0^+\),
\[
  \lambda_1^{(m)}(b) 
  = b - \frac{1}{2^{m - 1} (m - 1)!}b^{m + 1} 
  + \cO (b^{m + \frac32}).
\]
\end{theorem}

\begin{remark}
  The \(\mathcal{O}(b^{m + \frac{3}{2}})\) can likely be improved to
  \(\mathcal{O}(b^{m + 2})\).
\end{remark}

\begin{proof}[Proof of Theorem~\ref{thm:lambdamasymptotic}]
We organize our calculations to use the Temple inequality on \(\Hm - b\), to get
both the upper and lower bounds needed. We first work out the sizes of the terms
involved. Throughout the proof we use the shorthand notations
$\langle\cdot,\cdot\rangle$ and $\|\cdot\|$ for the inner product and the norm
in $L^2( (1,+\infty) ;r\dd r)$, respectively.  All order calculations
are done in the limit \(b \to 0^+\).

{\bfseries Step 1 (Quasimode).} Fix \(m\), and let \(\chi\) be the function defined as
\[
  \chi(r)=1+\frac{m-b/2}{m+b/2}r^{-2m}.
\]
We will define a quasimode \(\Psi\) in the domain of \(\Hm\),
\[
  \Psi(r) = \chi(r) f(r),\qquad f(r)=r^m\ee^{-br^2/4},
\]
with
\[
  \Psi'(1)=\bigl[ \chi'(1)+(m-b/2)\chi(1)\bigr]f(1)=0.
\]
We next show that
\begin{equation}\label{eq:Psinorm}
  \|\Psi\|^2 = \frac{2^m m!}{b^{m + 1}} + \cO (1/b).
\end{equation}
Since \(\Psi\) is real-valued, we can skip absolute values, and note that the
\(\Psi(r)^2\) can be written as a sum of the terms
\[
    A_1(r) = r^{2m}\ee^{-br^2/2},\quad
    A_2(r) = B_m^2r^{-2m}\ee^{-br^2/2},\quad
    A_3(r) = 2B_m\ee^{-br^2/2},
\]
where  $B_m=(m-b/2)/(m+b/2) = 1 + \cO(b)$. 

We will see that the integrals of \(A_2\) and \(A_3\) are small when compared to
the integral of $A_1$. We get
\[
  \int_1^{+\infty} A_1(r) r\dd r 
  = \int_1^{+\infty} r^{2m + 1}\ee^{-b r^2/2} \dd r 
  = \frac{2^m}{b^{m + 1}}\Gamma(m + 1,b/2)
  = \frac{2^m m!}{b^{m + 1}} + \cO (1).
\]
For the integral of $A_3$,  a change of variable yields
\[ 
  \int_1^{+\infty}A_3(r)r\dd r
  =
  \frac{2B_m}{b}\int_{b/2}^{+\infty}\ee^{-s}\dd s
  =
  \frac{2B_m\ee^{-b/2}}{b}
  =
  \cO(1/b),
\]
and for the integral of $A_2$,  we write by Hölder's inequality,
\[ 
  \int_1^{+\infty}A_2(r)r\dd r\leq B_m^2
  \Bigl(
    \int_1^{+\infty}r^{-4m} r\dd r
  \Bigr)^{1/2}
  \Bigl(
    \int_1^{+\infty}\ee^{-br^2} r\dd r
  \Bigr)^{1/2}
  =
  \frac{B_m^2\ee^{-b/2}}{\sqrt{4(2m-1)b}}.
\]
The claimed size of the norm in~\eqref{eq:Psinorm} therefore holds.

{\bfseries Step 2.} We next want to estimate the size of \(\norm{(\Hm -
b)\Psi}\). Observing that
\[
  \Hm f=bf,\quad f'=\Bigl(\frac{m}{r}-\frac{br}{2}\Bigr)f,
\]
we get
\[
  (\Hm - b)\Psi
  =
  -f\Bigl( \chi'' +\frac{1+2m}{r}\chi' - br\chi'\Bigr).
\] 
Now, we note that $\chi''+\frac{1+2m}{r}\chi'=0$,  hence
\begin{equation}\label{eq:Temple1}
  (\Hm - b)\Psi= -2b m B_m r^{-m}\ee^{-br^2/4},
\end{equation}
and, just as we estimated the integral of $A_2$ in the previous step, we get
eventually
\begin{equation}\label{eq:Temple2*}
  \norm{(\Hm - b)\Psi}=\cO \bigl(b^{\frac{3}{4}}\bigr).
\end{equation}

{\bfseries Step 3.} Now we estimate the inner product $\innerproduct{
(\Hm-b)\Psi,\Psi }$. Thanks to~\eqref{eq:Temple1}, we have
\[ 
  \innerproduct{ (\Hm-b)\Psi,\Psi} 
  = 
  -2b m B_m\int_1^{+\infty} \chi(r)\ee^{-br^2/2}r\dd r.
\]
We write the integral on the right-hand side as the sum of the following two
integrals
\[
  \begin{gathered}
    \int_1^{+\infty} \ee^{-br^2/2} r\dd r
    = 
    \ee^{-b/2}/b
     =
    \cO(1/b) \\
    B_m\int_1^{+\infty} r^{-2m}\ee^{-br^2/2} r\dd r
    =
    \cO \bigl(1/\sqrt{b}\bigr),
  \end{gathered}
\]
and consequently,  we end up with 
\begin{equation}\label{eq:Temple2}
\innerproduct{ (\Hm-b)\Psi,\Psi} = -2m + \cO \bigl(\sqrt{b}\bigr).
\end{equation}

{\bfseries Step 4 (Temple's inequality).} As $b\to 0^+$, we have a spectral gap,
thanks to Proposition~\ref{prop:spectrum-1-2} which yields
\[ 
  \lambda_2^{(m)}(b)-\lambda_1^{(m)}(b)\geq 2b+o(b).
\]
So we can apply Temple's inequality (see, e.g., \cite[Thm. 2]{H78}) for the
operator \(\Hm - b\),
\[
  \eta -\frac{\epsilon^2}{\beta-\eta}
  \leq
  \lambda_1^{(m)}(b) - b
  \leq
  \eta,
\]
where
\[
  \eta=\frac{\innerproduct{ (\Hm-b)\Psi,\Psi } }{\norm{\Psi}^2}
  \quad\text{and}\quad
  \epsilon^2=\frac{\norm{(\Hm-b)\Psi}^2}{\norm{\Psi}^2}-\eta^2 ,
\]
and $\beta=2b$.

Collecting the estimates in~\eqref{eq:Psinorm},~\eqref{eq:Temple1}
and~\eqref{eq:Temple2}, we get positive constants \(C_1,C_2\) and \(C_3\), as
well as \(\hat{b} > 0\) so that for \(0 < b < \hat{b}\), it holds that \(\beta
<\lambda_2^{(m)}(b)\),
\[
  \frac{b^{m+1}}{(m - 1)!2^{m - 1}} - C_1 b^{m + \frac{3}{2}}
  \leq
  \eta
  \leq
  \frac{b^{m+1}}{(m - 1)!2^{m - 1}} + C_2 b^{m + \frac{3}{2}}
\]
and \( 0\leq \epsilon^2 \leq C_3 b^{m + \frac{5}{2}} \). The promised asymptotic
formula follows.
\end{proof}

\subsection{Ordering of the eigenvalues}

Our next aim is to prove the existence of \(b_0 > 0\) such that, for all \(0 < b
< b_0\), the eigenvalues are ordered as
\[
  \lambda_1^{(1)}(b) < \lambda_1^{(2)}(b) < \lambda_1^{(3)}(b) < \ldots < b.
\]
To achieve this, we will use the asymptotic expansion of each
\(\lambda_1^{(m)}(b)\) as \(b \to 0^+\), as well as some information about
intersections.

\subsubsection{Intersection of band functions}
With the asymptotic formula in Theorem~\ref{thm:lambdamasymptotic} at hand, we
can now with the help of Proposition~\ref{prop:intersections} find that the
first eigenvalues \(\lambda_1^{(m)}(b)\) do not intersect for small \(b\).

First we note that the asymptotic formula gives the existence of constants
\(\hat{b}_m\) such that if \(0 < b < \hat{b}_m\) then
\[
  \lambda_1^{(m - 1)}(b) < \lambda_1^{(m)}(b).
\]
The next proposition states an explicit bound for \(\hat{b}_m\).

\begin{proposition}\label{prop:ordering}
  Assume that \(m \geq 2\) and \(0 < b \leq 2m + 1 - \sqrt{8m + 1}\). Then
  \[
    \lambda_1^{(m - 1)}(b) < \lambda_1^{(m)}(b).
  \]
\end{proposition}
\begin{proof}
  From Proposition~\ref{prop:intersections} we know that if \(\lambda_1^{(m -
  1)}(b) = \lambda_1^{(m)}(b) = \lambda < b\) then
  \[
    (b/2 - m)\bigl(b/2 - (m - 1)\bigr) < b.
  \]
  We recall that both eigenvalues \(\lambda_1^{(m - 1)}(b)\) and
  \(\lambda_1^{(m)}(b)\) are smaller than \(b\) only if \(b < 2(m - 1)\). For
  these values of \(b\) the displayed inequality above is equivalent to \(b > 2m
  + 1 - \sqrt{8m + 1}\). Hence, if \(b \leq 2m + 1 - \sqrt{8m + 1}\) then
  \(\lambda_1^{(m - 1)}(b)\) and \(\lambda_1^{(m)}(b)\) cannot be equal, and
  since by Theorem~\ref{thm:lambdamasymptotic} it holds that  \(\lambda_1^{(m -
  1)}(b) < \lambda_1^{(m)}(b)\) for small \(b\), this must still be true for all
  \(b\) with \(0 < b \leq 2m + 1 - \sqrt{8m + 1}\).
\end{proof}

\begin{corollary}\label{cor:ordering}
  If \(0 < b < 5 - \sqrt{17}\) then
  \[
    \lambda_1^{(1)}(b) < \lambda_1^{(2)}(b) < \lambda_1^{(3)}(b) < \ldots < b.
  \]
\end{corollary}

\begin{proof}
This is immediate from  Proposition~\ref{prop:ordering} and the fact that \(m
\mapsto 2m + 1 - \sqrt{8m + 1}\) is increasing for \(m \geq 2\).
\end{proof}

\subsubsection{Uniform convergence of eigenvalues}

\begin{proposition}
For \(0 < b < 5 - \sqrt{17}\),  it holds that \(\lambda_1^{(m)}(b) \to b\)
uniformly as \(m \to +\infty\).
\end{proposition}

\begin{proof}
We know that the eigenvalues obey in this interval
\[
  \lambda_1^{(1)}(b) < \lambda_1^{(2)}(b) < \lambda_1^{(3)}(b) < \ldots < b.  
\]
The eigenvalues must converge to \(b\) pointwise as \(m \to +\infty\) since
\(b\) is the bottom of the essential spectrum of \(\Hb\). The convergence is
uniform by the Dini theorem.
\end{proof}

\begin{remark}
The convergence will in fact be uniform on any interval \((0,\hat{b})\),
\(\hat{b} > 0\), since we have a good control on where the intersections occur.
\end{remark}

\subsection{Proof of Theorem~\ref{thm:disk}} 
Assume that $R=1$. By (2) in Proposition~\ref{prop:disk-m}, the ground state is
non-radial. Moreover, by (1) and (4) in Proposition~\ref{prop:disk-m}, and by
Corollary~\ref{cor:ordering}, we get for \(0 < b < 5 - \sqrt{17}\) that
\[
  \lambda_1(b,\disk_1^{\ext})
  =
  \lambda_1^{(1)}(b)
  <
  \lambda_2(b,\disk_1^{\ext})
  =
  \lambda_2^{(2)}(b)
  <
  \ldots
  <
  b.
\]
Finally, the announced two term asymptotics of $\lambda_k(b,\disk_1^{\ext})$
follows from Theorem~\ref{thm:lambdamasymptotic}.

\section{Optimality of  the disk's exterior}\label{sec:loc-opt}
The purpose of this section to prove Theorems \ref{thm:ub-general} and \ref{thm:loc-opt}. Recall that
$\disk_1\subset\reals^2$ denotes the unit disk centred at the origin and that
$\disk_1^{\ext} := \reals^2\setminus\overline{\disk_1}$ is the exterior of
$\disk_1$.

Throughout this section,  we suppose that  $0 < b < 5-\sqrt{17}$. Thanks to
Corollary~\ref{cor:ordering}, and Proposition~\ref{prop:disk-m}, we know that
\begin{equation}\label{eq:ev-loc-op}
  \lambda_1(b,\disk_1^{\ext})
  =
  \lambda_1^{(1)}(b),
\end{equation}
is the lowest eigenvalue of the fibre operator $\cH^{(1)}(b)$ introduced in
\eqref{eq:Hm-b}. 
In particular, we know that, as \(b \to 0^+\),
\begin{equation}
  \lambda_1(b,\disk_{1}^{\ext})
  =
  b - b^2 + o(b^2).
\end{equation}
This is what we need to compare with for general domains, when addressing
optimality of the exterior of the unit disk in the limit $b\to 0^+$.

\subsection{Trial function}

The ground state of $\cH$ on $\disk_1^{\ext}$ can be represented in polar
coordinates as 
\begin{equation}\label{eq:def-u-polar}
  u(r,\theta) := f(r)\ee^{\ii\theta}, 
\end{equation}
where $f$ is the positive normalized  eigenfunction of the fibre operator
$\cH^{(1)}(b)$ corresponding to its lowest eigenvalue $\lambda :=
\lambda_1^{(1)}(b) < b$. In particular, we know that $f\in L^2((1,+\infty);r\dd
r)$ satisfies the ordinary differential equation 
\begin{equation}\label{eq:differential}
	- f''(r) 
  - \frac{1}{r}f'(r) 
  + \left(\frac{1}{r}-\frac{br}{2}\right)^2f(r) 
  = \lambda f(r),
	\qquad r > 1,
\end{equation}
and the boundary condition $f'(1) = 0$. By standard properties of ordinary
differential equations, we can construct a $C^\infty$-smooth extension $f$ to
$(0,\infty)$, which satisfies the same differential
equation~\eqref{eq:differential} for all $r > 0$.  This function can be expressed
in terms of a confluent hypergeometric function (Whittaker function), see
Appendix~\ref{sec:app}. Without any danger of confusion we will denote the
respective extension again by $f$. Let us also introduce the extension of $u$ to
$\reals^2\setminus\{0\}$,
\begin{equation}\label{eq:def-trial-fct}
	v(r,\theta) = f(r)\ee^{\ii\theta},\qquad r > 0,\ \theta \in [0,2\pi).
\end{equation}
It is not hard to see that
\begin{equation}\label{eq:v_property}
	\bigl((-\ii \nabla - b{\Ab})^2 v\bigr)(x) = \lambda v(x),\qquad \text{for all } x\in\reals^2\setminus\{0\}.
\end{equation}

\subsection{Upper bound for general exterior domains}

Using the trial function in~\eqref{eq:def-trial-fct},  we can prove an upper
bound on the lowest eigenvalue for an exterior domain not containing the origin
of $\reals^2$.

\begin{proposition}\label{prop:ub-general} 
  Suppose that $\Omegaext=\reals^2\setminus\overline{\Omega}$, where $\Omega$ a
  bounded simply connected domain with a $C^\infty$-smooth boundary and such that
  $0\in\Omega$.  If $b < 5-\sqrt{17}$, then the lowest eigenvalue in
  $\Omegaext$ satisfies
  \begin{equation}\label{eq:minmax}
    \lambda_1(b,\Omegaext)
    \leq
    \lambda_1(b,\disk_1^{\ext})
    -\frac{\displaystyle\int_{\partial\Omega}\big(\nu\cdot \nabla v\big) \overline{v}\dd \sigma}
          {\displaystyle\int_{\Omegaext }|v|^2\dd x},
  \end{equation}
	where $v$ is the function introduced in~\eqref{eq:def-trial-fct}, and $\nu$
	is the outer unit normal for $\Omega$. 
\end{proposition}

\begin{proof}
By construction, we have  $v|_{\Omegaext }\in H^1_{{\Ab}}(\Omegaext )$. For the
sake of brevity we will denote the restriction of $v$ to $\Omegaext $ again by
$v$ as no confusion can arise. Using the min-max principle and integrating by
parts we get
\[
  \begin{aligned}
    \lambda_1(b,\Omegaext ) 
    &\le 
    \frac{\displaystyle\int_{\Omegaext }|(-\ii\nabla - b{\Ab})v|^2\dd x}
        {\displaystyle\int_{\Omegaext }|v|^2\dd x}\\
    &= 
    \frac{\displaystyle\int_{\Omegaext }(-\ii\nabla - b{\Ab})^2v \overline {v}\dd x}
        {\displaystyle\int_{\Omegaext }|v|^2\dd x}
    -
    \frac{\displaystyle\int_{\partial\Omega}\big(\nu\cdot(\nabla -\ii b{\Ab})v\big) \overline{v}\dd \sigma}
        {\displaystyle\int_{\Omegaext }|v|^2\dd x}\\
    &= 
    \lambda_1(b,\disk_1^{\ext})
    - 
    \frac{\displaystyle\int_{\partial\Omega}\big(\nu\cdot \nabla v\big) \overline{v}\dd \sigma}
        {\displaystyle\int_{\Omegaext }|v|^2\dd x},		
  \end{aligned}	
\]
where we used in the last step that $v$ satisfies the differential equation
in~\eqref{eq:v_property} for the first term, whereas the magnetic contribution
in the boundary term vanishes. Indeed, integrating by parts we get
\[
  \begin{aligned}
    \int_{\partial\Omega} (\nu\cdot \Ab)|v|^2\dd \sigma
    & = 
    -\int_{\Omegaext }\Div (\Ab |v|^2)\dd x\\
    & =
    -\int_{\Omegaext }(\Div \Ab)|v|^2\dd x
    -\int_{\Omegaext }\Ab \cdot \nabla(|v|^2)\dd x =0,
  \end{aligned}
\]  
where in the last step we used that $\Div \Ab = 0$ and that $\Ab(x)$ is
orthogonal to $\nabla(|v|^2)(x)$ for any $x\in\Omegaext$,  since \(|v|\) is
radial.
\end{proof}

\subsection{Star-shaped domains}

When $\Omega$ is star-shaped and contains the origin,  we can express the
boundary term in~\eqref{eq:minmax} in terms of the following function
\begin{equation}\label{eq:def-F-a}
	F_{\alpha}(r) 
  \coloneqq 
  r^{\alpha} f'(r^{\alpha})f(r^{\alpha}), \quad (\alpha\neq 0,~r > 0),
\end{equation}
where $f$ is the function appearing in~\eqref{eq:def-trial-fct}. Notice that
$F_1(r)=F_{\alpha}(r^{1/\alpha})$, for any $\alpha\neq 0$.

\begin{proposition}\label{prop:bnd-term} 
  Suppose that $\Omega$ is defined by polar coordinates as $\Omega =
  \{(r,\theta)\colon r < \rho(\theta)\}$, where $\rho\colon \reals/[0,2\pi]\to
  \reals_+$ is a $C^\infty$-smooth function. Then,
  \begin{equation}\label{eq:bterm}
    \int_{\partial\Omega} (\nu\cdot \nabla v)\overline{v}\dd\sigma
    = \int_0^{2\pi} F_1\bigl(\rho(\theta)\bigr) \dd\theta,
  \end{equation}
  where $v$ is the function introduced in~\eqref{eq:def-trial-fct}, and $F_1$
  is the function introduced in~\eqref{eq:def-F-a}.
\end{proposition}

\begin{proof}
We can write the outer unit normal vector field to $\partial\Omega$  as a function of the angular variable (cf.~\cite[Eq.~(3.10)]{KL23})
\begin{equation}\label{eq:normal}
	\nu(\theta) = \frac{1}{\sqrt{\rho^2(\theta)+ [\dot\rho(\theta)]^2}}
	\begin{pmatrix}
    \rho(\theta)\cos\theta +\dot{\rho}(\theta)\sin\theta,\rho(\theta)\sin\theta - \dot{\rho}(\theta)\cos\theta
  \end{pmatrix}^\top.
\end{equation}
The gradient of the trial function $v$ in polar coordinates can be written as
\begin{equation}\label{eq:gradtrial}
	(\nabla v)(r,\theta) = f'(r) \ee^{\ii \theta}\begin{pmatrix}\cos\theta\\
		\sin\theta\end{pmatrix} +  \frac{\ii f(r)\ee^{\ii\theta}}{r}
	\begin{pmatrix}
    -\sin\theta\\
		\cos\theta
  \end{pmatrix}.
\end{equation}
Using~\eqref{eq:normal} and~\eqref{eq:gradtrial} we can express the boundary
term as 
\[
\begin{split}
	\int_{\partial\Omega} (\nu\cdot \nabla v)\overline{v}\dd\sigma &=
	\int_0^{2\pi}\left(\rho(\theta) f'(\rho(\theta)) -\ii \frac{\dot\rho(\theta)}{\rho(\theta)}f(\rho(\theta))\right) f(\rho(\theta))\dd \theta\\
	& = \int_0^{2\pi} \rho(\theta)f'(\rho(\theta))f(\rho(\theta))\dd\theta.\qedhere
\end{split}
\]
\end{proof}

Thanks to Propositions~\ref{prop:ub-general} and~\ref{prop:bnd-term},  we need
to determine the sign of the integral involving $F_1$  in order to finish the
proof of Theorem~\ref{thm:loc-opt}.  Towards that aim,  the  next proposition
will be handy.
\begin{proposition}\label{prop:F-a}
If $\alpha<0$,  then the function $F_\alpha$ introduced in~\eqref{eq:def-F-a}
satisfies
\[
	F_{\alpha}''(1) \ge |f(1)|^2\left[2\alpha^2\left(
	\frac{b^2}{2}-2b\right) -\alpha\left(1-2b +\frac{b^2}{4}\right)\right].
\]
In particular,  if 
\[
	0 <b < 4-2\sqrt{3}\quad \text{ and }\quad\frac{1-2b+\frac{b^2}{4}}{b^2-4b} < \alpha < 0,
\]
then $F_{\alpha}''(1) > 0$.
\end{proposition}

\begin{proof}
Differentiating \(F_{\alpha}\) twice, inserting \(r = 1\), and using \(f'(1) =
0\), we get
\begin{equation}\label{eq:F2}
	F_{\alpha}''(1) = (3\alpha^2-\alpha) f''(1)f(1) + \alpha^2f'''(1)f(1).
\end{equation}
From the differential equation~\eqref{eq:differential} employing that $f'(1) =
0$ we derive
\begin{equation}\label{eq:f1}
	f''(1) = \left[\left(1-\frac{b}{2}\right)^2-\lambda\right] f(1).
\end{equation}
Differentiating~\eqref{eq:differential} with respect to $r$ 
and inserting \(r = 1\), we get
\begin{equation}\label{eq:f2}
  f'''(1) + f''(1) = -2\left(1-\frac{b^2}{4}\right)f(1).
\end{equation}
Substituting~\eqref{eq:f1} and \eqref{eq:f2} into~\eqref{eq:F2}
\[
\begin{aligned}
	F_{\alpha}''(1) 
  & =
  |f(1)|^2
  \left[
    (2\alpha^2-\alpha)\left(1-b+\frac{b^2}{4}-\lambda\right) 
    - 2\alpha^2\left(1 - \frac{b^2}{4}\right)
  \right].
\end{aligned}	
\]
For $\alpha < 0$ we get using the inequality $\lambda < b$,
\[
	F_{\alpha}''(1) \ge |f(1)|^2\left[2\alpha^2\left(
	\frac{b^2}{2}-2b\right) -\alpha\left(1-2b +\frac{b^2}{4}\right)\right].
\]
Writing
\[
  2\alpha^2
  \left(
	  \frac{b^2}{2} - 2b
  \right) 
  - \alpha 
  \left(
    1-2b +\frac{b^2}{4}
  \right)
  =
  -\alpha
  \left[
    1-2b+\frac{b^2}{4}-2\alpha\Bigl(\frac{b^2}{2}-2b\Bigr)
	\right]
\]
then the conditions $0 <b < 4-2\sqrt{3}$ and $\frac{1-2b+\frac{b^2}{4}}{b^2-4b}
< \alpha < 0$ yield 
\[
  1-2b +\frac{b^2}{4} > 0
  \quad\text{ and }\quad 
  1-2b+\frac{b^2}{4}-2\alpha\Bigl(\frac{b^2}{2}-2b\Bigr) > 0,
\]
and consequently $F_{\alpha}''(1)>0$.
\end{proof}

\subsection{Proof of Theorem~\ref{thm:loc-opt}}

We will prove the theorem for the unit radius of the disk ($R = 1$). The general
case follows by scaling. The star-shaped domain $\Omega$ is  parametrized by a
smooth mapping $\reals/[0,2\pi]\ni\theta\mapsto\rho(\theta)$ as $\Omega =
\{(r,\theta)\colon r <\rho(\theta)\}$,  and in this case   the $p$-moment of
$\Omega$ can be computed as follows:
\begin{equation}\label{eq:moment}
	\cI_p(\Omegaext ) 
  = 
  \int_0^{2\pi}\int_{\rho(\theta)}^\infty r^{1-p}\dd r 
  = 
  \frac{1}{p-2}\int_0^{2\pi} (\rho(\theta))^{2-p}\dd\theta.
\end{equation}
Recall that we impose the constraint $\cI_p(\Omegaext)=\cI_{p}(\disk_1^{\ext})$.  
Let us  set
\[
	\alpha = \frac{1}{2-p}.
\]
By the assumptions on $p$ and $b$ in the formulation of
Theorem~\ref{thm:loc-opt}, the conditions in Proposition~\ref{prop:F-a} are
fulfilled, and we have $F''_{\alpha}(1)>0$.  By smoothness of $F_{\alpha}$ there
exists $\varepsilon\in (0,1)$ such that $F_{\alpha}''(r) > 0$ for all $r\in
[1-\varepsilon,1+\varepsilon]$,  hence $F_{\alpha}$ is convex on
$[1-\varepsilon,1+\varepsilon]$.

Thus, there exists $r_\star \in (0,1)$ such that $\rho^{1/\alpha}\in
(1-\varepsilon,1+\varepsilon)$ provided that $1-r_\star \le \rho \le 1+
r_\star$. Using Jensen's inequality and~\eqref{eq:moment} we get
\[
\begin{aligned}
	\int_0^{2\pi}F_1(\rho(\theta))\dd \theta &= 
	\int_0^{2\pi}
	F_{\alpha}((\rho(\theta))^{1/\alpha})\dd \theta\\
	&
  \ge 2\pi F_{\alpha}\left(\frac{1}{2\pi}\int_0^{2\pi}
	(\rho(\theta))^{1/\alpha}\dd \theta\right)\\
	&
	=2\pi F_{\alpha}\left(\frac{p-2}{2\pi}\cI_p(\Omegaext )\right)\\
	&
  = 2\pi F_{\alpha}\left(\frac{p-2}{2\pi}\cI_p(\disk_1^{\ext})\right)\\	
	&
  = 2\pi F_{\alpha}(1) = 0.
\end{aligned}	
\]
Hence, the boundary term in~\eqref{eq:minmax} is positive by
Proposition~\ref{prop:bnd-term},  and the desired isoperimetric inequality
follows from Proposition~\ref{prop:ub-general}. This finishes the proof of
Theorem~\ref{thm:loc-opt} for $R = 1$.   The general case follows by scaling.

\begin{remark}
	The fixed area constraint corresponds to the choice $\alpha = \frac12$ in the
	above proof of Theorem~\ref{thm:loc-opt} (due to the formula $|\Omega| =
	\frac12\int_0^{2\pi}\rho^2(\theta)\dd \theta$ for star-shaped $\Omega$). In
	this case, we get $F_{1/2}''(1) = -\frac12|f(1)|^2 (1 - b^2/4)$ is positive
	provided that $b > 2$. However, for such magnetic fields the ground state of
	the magnetic Neumann Laplacian on the exterior of the unit disk does not
	correspond to the fibre labelled by $m = 1$,  and we no more
	have~\eqref{eq:ev-loc-op}.   Thus, our method is not directly applicable.  
\end{remark}

\subsection{Proof of Theorem~\ref{thm:ub-general}}
We can derive an upper bound when $b\to 0^+$ and $\Omega=\{(r,\theta)\colon
r < \rho(\theta)\}$ is star-shaped and contains $0$, by modifying the trial
state in~\eqref{eq:def-trial-fct} as follows
\begin{equation}\label{eq:trial-fct}
  u(r,\theta) 
  = 
  \Psi(r)\ee^{\ii\theta},
  \qquad 
  \Psi(r) = \big(r+cr^{-1} \big)\ee^{-br^2/4},
\end{equation}
with some parameter $c\in\mathbb{R}_+$ to be determined at a later stage.
The function $\Psi$ is reminiscent of the one appearing in the proof of
Theorem~\ref{thm:lambdamasymptotic} for $m=1$.  It satisfies the differential
equation
\[
  - \Psi''
  - \frac1r \Psi'
  + \Bigl(\frac{1}{r}-\frac{br}{2}\Bigr)^2\Psi
  - b\Psi
  =
  - 2b c r^{-1}\ee^{-br^2/4}
  \leq 0.
\]
Fixing a positive constant $\varepsilon$,  similar calculations as in the proof
of Theorem~\ref{thm:lambdamasymptotic} yield the following asymptotics as $b \to
0^+$,
\[ 
  \begin{gathered}
    \int_\varepsilon^{+\infty}|\Psi|^2r\dd r
    =
    \frac{2}{b^2}+\cO(1/b),\\
    \int_\varepsilon^{+\infty}
      \Bigl[
        - \Psi'' 
        - \frac{1}{r}\Psi' 
        + \Bigl(\frac1r-\frac{br}2\Bigr)^2 \Psi
        - b\Psi
      \Bigr]
      \Psi \,r\dd r
    =
    - 2c+ \cO(\sqrt{b}).
  \end{gathered}
\]
Observing that there are positive constants $r_1<r_2$ such that
\[
  \disk_{r_1}\subset\Omega\subset \disk_{r_2},
\]
we get by a straightforward comparison with the integrals over
$\disk_{r_1}^{\ext}$ and $\disk_{r_2}^{\ext}$,
\[ 
  \|u\|_{L^2(\Omegaext)}^2=\frac{4\pi}{b^2}+\cO(1/b),
  \quad 
  \innerproduct[\big]{ [(-\ii\nabla-b\Ab)^2-b]u,u }_{L^2(\Omegaext)}
  =-4\pi c+\cO(\sqrt{b}).
\]
Arguing as in Propositions~\ref{prop:ub-general} and \ref{prop:bnd-term}, we get
by the min-max principle
\[
  \begin{aligned}
  \lambda_1(b,\Omegaext)
  &
  \leq 
  \frac{\innerproduct{ (-\ii\nabla-b\Ab)^2u,u }_{L^2(\Omegaext)}}
       {\|u\|_{L^2(\Omegaext)}^2}
  - \frac{\displaystyle\int_{\partial\Omega}\nu\cdot\nabla u\,\overline{u}\dd\sigma}
         {\|u\|_{L^2(\Omegaext)}^2}\\
  &
  =
  b - b^2c - 
  \frac{b^2}{4\pi}
  \int_0^{2\pi}\rho(\theta)\Psi'(\rho(\theta))\Psi(\rho(\theta))\dd \theta
  +\cO(b^{5/2}).
  \end{aligned}
\]
A routine calculation shows that,  as $b\to 0^+$,  
\[
  r\Psi'(r)\Psi(r)=r^2-c^2r^{-2}+\cO(b)
\]
uniformly with respect to $r$,  as long as  $r$ is in a compact interval of
$\reals_+$.  Consequently,
\[ 
  \int_0^{2\pi}\rho(\theta)\Psi'(\rho(\theta))\Psi(\rho(\theta))\dd\theta
  =
	 2\bigl(|\Omega|-c^2\cI_4(\Omegaext)\bigr)
  + \cO(b).
\] 
We therefore have proved that
\begin{equation}\label{eq:ub-gen-dom}
  \lambda_1(b,\Omegaext)
  \leq 
  b - \biggl(
        c + \frac{|\Omega|-c^2\cI_4(\Omegaext)}{2\pi} 
      \biggr)b^2 + \cO(b^{5/2}).
\end{equation}
Optimizing the bound with respect to the parameter $c$ we easily get that the optimal value is $c = \frac{\pi}{\cI_4(\Omegaext)}$
and we obtain the following:
\begin{equation}\label{eq:proof.ub-general}
	\lambda_1(b,\Omegaext)\leq b-\cG(\Omega)b^2+\cO(b^{5/2}),
\end{equation}
where
\[
\cG(\Omega)=\frac{|\Omega|}{2\pi}+\frac{\pi}{2\cI_4(\Omegaext)}.
\]
This bound is consistent with the conclusion of Theorem~\ref{thm:disk},  since  $\cG(\Omega)=R^2$ if $\Omega=\disk_R$ is the disk of radius $R$.

In
the remaining part of the proof we assume that $\Omega$ is not a disk centred
at the origin. By Hölder's inequality,
\[
2\pi
<
\Bigl(\int_0^{2\pi}\rho^2\dd\theta\Bigr)^{1/2}
\Bigl(\int_0^{2\pi}\rho^{-2}\dd\theta\Bigr)^{1/2},
\]
which reads as $\pi^2 < |\Omega|\cI_4(\Omegaext)$, 
Thus, we always have
\[
\frac{\pi}{\cI_4(\Omegaext)}
<
\cG(\Omega)
<
\frac{|\Omega|}{\pi}.
\]
In particular,  if we fix $R_\star$ such that
$\cI_4(\disk_{R_\star}^\ext)=\cI_4(\Omegaext)$,  we have
\[
\cG(\Omega)
>
\cG(\disk_{R_\star}),
\]
and we get from \eqref{eq:proof.ub-general} and Theorem~\ref{thm:disk} that
$\lambda_1(b,\Omegaext) < \lambda_1(b,\disk_{R_\star}^{\ext})$, for $b$
sufficiently small,  by which the proof of Theorem~\ref{thm:ub-general} is complete.

\begin{remark}
However, if we fix the area by choosing $R$ so that $|\disk_R|=|\Omega|$,
then $\cG(\Omega) < \cG(\disk_R)$ and \eqref{eq:proof.ub-general} is not
sufficient to compare $\lambda_1(b,\Omegaext)$ and
$\lambda_1(b,\disk_{R}^{\ext})$.
\end{remark}

\subsection*{Acknowledgments} {\small
The first listed author (A.~K.) is partially supported by The Chinese University
of Hong Kong, Shenzhen (grant UDF01003322). The second listed author (V.~L.) is
supported by the  Czech Science Foundation (GA\v{C}R) within the project
21-07129S.

We initiated this work when the first listed author  visited Lund University and
the Nuclear Physics Institute, Czech Academy of Sciences, \v{R}e\v{z}, Czech
Republic. We acknowledge financial support through Knut and Alice Wallenberg
Foundation (grant KAW 2021.0259).}

\appendix

\section{Intersection of eigenvalues}\label{sec:app}

\subsection{Introduction}

We consider here the self-adjoint realization of the operator
\begin{equation}\label{eq:mdiffekv}
  \Hm(b) 
  = 
  - \frac{\dd^2}{\dd r^2}
  - \frac{1}{r}\frac{\dd}{\dd r} 
  + \Bigl(\frac{m}{r} - \frac{br}{2}\Bigr)^2,
\end{equation}
either in \(L^2((0,1),r\dd r)\) or \(L^2((1,+\infty),r\dd r)\), with Neumann
condition at \(r = 1\). We denote in both cases by \(\lambda_1^{(m)}(b)\) the
lowest eigenvalue of $\Hm(b)$. 

\subsection{The interior of the disk}
We start by working in the interior of the disk. The following result is old,
but not so well-known, since it is in the original paper only mentioned as a
comment, without a proof.

\begin{theorem}[{\cite{sj}}]
  Assume that \(b > 0\) and \(m \geq 1\). If \(\lambda_{1}^{(m - 1)}(b) =
  \lambda_1^{(m)}(b) < b\) then \(\lambda_1^{(m)}(b) = (b/2-m)[b/2-(m-1)]\).
\end{theorem}

The condition that the eigenvalues should be less than $b$ is merely to enforce
the eigenvalues to be the bottom ones, and to be sure that they are not
eigenvalues of the Dirichlet operator.

\begin{proof}
The result relies on several identities for confluent hypergeometric functions,
we will use the so-called Whittaker functions.

The Whittaker differential equation reads
\begin{equation}
  \frac{\dd^2y}{\dd z^2}
  + \Bigl(-\frac{1}{4}+\frac{\kappa}{z}+\frac{1/4-\mu^2}{z^2}\Bigr)y 
  = 0.
\end{equation}
It has two linearly independent solutions, denoted by \(\M(\kappa,\mu,z)\)
(well-defined if \(2\mu\) is not a negative integer, which it will not be in our
case) and \(\W(\kappa,\mu,z)\). According to~\cite[(13.14.14)]{dlmf},
\begin{equation}\label{eq:Mbounded}
  \M(\kappa,\mu,z) = z^{\mu+1/2} (1 + O(z)) \qquad (z\to 0).
\end{equation}
With a change of variables $z = br^2/2$ and with a scaling by $r$ of the
independent function, one finds that the general solution to $\Hm(b)u = \lambda
u$ is given by
\begin{equation}
  u(r) = c_1\frac{1}{r} \M\parentes*{\Lambda,\frac{m}{2},\frac{b r^2}{2}}
       + c_2\frac{1}{r} \W\parentes*{\Lambda,\frac{m}{2},\frac{b r^2}{2}}
\end{equation}
For notational convenience in the continuation, we have introduced the constant
\[
  \Lambda=\frac{m b + \lambda}{2b}.
\]
Since $\M$ is bounded at $0$ according to~\eqref{eq:Mbounded}, and since the
differential equation has a regular singularity at the origin, $\W$ becomes
singular as $r \to 0^+$. We are forced to take $c_2 = 0$. Let us also set $c_1 =
1$, and consider below the eigenfunction
\begin{equation}\label{eq:eigenfunction}
  u(r)=\frac{1}{r} \M\Bigl(\Lambda,\frac{m}{2},\frac{br^2}{2}\Bigr).
\end{equation}
The Whittaker function $\M$ satisfies~\cite[(13.15.15)]{dlmf}
\begin{equation}\label{eq:diffM}
  \frac{d}{dz}\M(\kappa,\mu,z)
  =
  \parentes*{\frac{1}{2}-\frac{\kappa}{z}}\M(\kappa,\mu,z)
  + \frac{1+2\kappa+2\mu}{2z}\M(\kappa+1,\mu,z).
\end{equation}
Implementing this, we find that $u'(1)=0$ precisely when
\begin{equation}\label{eq:form}
  \bigl(b^2 - 2b(m+1)-2\lambda\bigr)
  \M\parentes*{\Lambda,\frac{m}{2},\frac{b}{2}}
  + 2\bigl(b + 2bm + \lambda\bigr)
  \M\parentes*{\Lambda + 1,\frac{m}{2},\frac{b}{2}} 
  = 0.
\end{equation}
Replacing $m$ by $m-1$ we get
\begin{equation}\label{eq:formm1}
  \begin{multlined}
    \bigl(b^2 - 2bm - 2 \lambda\bigr)
    \M\parentes*{\Lambda - \frac{1}{2},\frac{m - 1}{2}, \frac{b}{2}}\\
    + 2\bigl(-b + 2bm + \lambda\bigr)
    \M\parentes*{\Lambda + \frac{1}{2},\frac{m - 1}{2}, \frac{b}{2}} 
    = 0.
  \end{multlined}
\end{equation}
We note that the third arguments in the Whittaker functions in~\eqref{eq:form}
and~\eqref{eq:formm1} are the same. We will \enquote{massage} the first
equation~\eqref{eq:form} in order to have Whittaker functions with the same
arguments in both equations.

In the first step we will make sure that the second arguments become the same.
To do this, we will use formula~\cite[(13.15.4)]{dlmf}, which says that
\begin{equation}\label{eq:Mfinal}
  \M(\kappa,\mu,z) = \frac{2\mu}{\sqrt{z}}\parentes*{\M\parentes{\kappa-\frac{1}{2},\mu-\frac{1}{2},z}-\M\parentes{\kappa+\frac{1}{2},\mu-\frac{1}{2},z}}.
\end{equation}
Implementing this in~\eqref{eq:form}, the equation transforms into (here we have
skipped the non-zero factor $m/\sqrt{b/2}$ from both terms)
\begin{equation}\label{eq:nym}
  \begin{multlined} 
    \bigl(b^2 - 2b(m+1) - 2\lambda\bigr)
    \bracket*{\M\parentes*{\Lambda-\frac{1}{2},\frac{m-1}{2},\frac{b}{2}}
    - \M\parentes*{\Lambda+\frac{1}{2},\frac{m-1}{2},\frac{b}{2}}}\\
    \quad +2\bigl(b + 2bm + \lambda\bigr)
    \bracket*{\M\parentes*{\Lambda+\frac{1}{2}, \frac{m - 1}{2}, \frac{b}{2}}
    - \M\parentes*{\Lambda+\frac{3}{2}, \frac{m - 1}{2}, \frac{b}{2}}} = 0.
  \end{multlined}
\end{equation}
We need a formula that makes the first argument comparable as well. To this aim,
we will use the formula \cite[(13.15.1)]{dlmf},
\begin{equation}\label{eq:Mfirstarg}
    \parentes*{\kappa + \mu + \frac{1}{2}}\M(\kappa + 1, \mu, z)
  = \parentes*{\mu - \kappa + \frac{1}{2}}\M(\kappa - 1, \mu, z)
    +(2\kappa - z)\M(\kappa, \mu, z).
\end{equation}
We will soon implement it in~\eqref{eq:nym}, with
\[
  \kappa = \Lambda + \frac{1}{2},
  \quad
  \mu = \frac{m - 1}{2},
  \quad
  z = \frac{b}{2}.
\]
This gives
\[
  \mu + \kappa + \frac{1}{2} = \Lambda + \frac{m + 1}{2},\quad
  \mu - \kappa + \frac{1}{2} = \frac{m - 1}{2} - \Lambda,\quad
  2\kappa - z = 2\Lambda + 1 - \frac{b}{2},
\]
and 
\[
  \begin{aligned}
    \parentes*{\Lambda + \frac{m + 1}{2}}
    \M\parentes*{\Lambda + \frac{3}{2}, \frac{m - 1}{2}, \frac{b}{2}}
    & =
    \parentes*{\frac{m - 1}{2} - \Lambda}
    \M\parentes*{\Lambda - \frac{1}{2}, \frac{m - 1}{2}, \frac{b}{2}} \\
    &\quad+
    \parentes*{2\Lambda + 1 - \frac{b}{2}}
    \M\parentes*{\Lambda + \frac{1}{2}, \frac{m - 1}{2}, \frac{b}{2}}.
  \end{aligned}
\]
We are now ready to insert this into~\eqref{eq:nym},
\begin{equation}\label{eq:nyarem}
  \begin{multlined}
    \bigl(b^2 - 2b(m + 1) - 2\lambda\bigr)\parentes*{\Lambda + \frac{m + 1}{2}}\\
    \times\bracket*{\M\parentes{\Lambda - \frac{1}{2}, \frac{m - 1}{2}, \frac{b}{2}}
    - \M\parentes*{\Lambda + \frac{1}{2}, \frac{m - 1}{2}, \frac{b}{2}}} \\
     + 2\bigl(b + 2bm + \lambda\bigr)
    \parentes*{\Lambda + \frac{m + 1}{2}}
    \M\parentes*{\Lambda + \frac{1}{2}, \frac{m - 1}{2}, \frac{b}{2}} \\
    -2\bigl(b + 2bm + \lambda\bigr)
    \parentes*{\frac{m - 1}{2} - \Lambda}
    \M\parentes*{\Lambda - \frac{1}{2}, \frac{m - 1}{2}, \frac{b}{2}} \\
    -2\bigl(b + 2bm + \lambda\bigr)
    \parentes*{2\Lambda + 1 - \frac{b}{2}}
    \M\parentes*{\Lambda + \frac{1}{2}, \frac{m - 1}{2}, \frac{b}{2}} = 0.
  \end{multlined}
\end{equation}
Collecting the two different Whittaker expressions, and cancelling the non-zero
common factor $(b + 2bm + \lambda)/2$, we find that
\begin{equation}\label{eq:finalm1}
   (b - 2m)\M\parentes*{\Lambda - \frac{1}{2}, \frac{m - 1}{2}, \frac{b}{2}}
  +(b + 2m)\M\parentes*{\Lambda + \frac{1}{2}, \frac{m - 1}{2}, \frac{b}{2}} = 0.
\end{equation}
Now assume that we have eigenvalues $\lambda_1^{(m)}(b)=\lambda_1^{(m - 1)}(b) =
\lambda$. This means that we have
\[
  \begin{bmatrix}
    b^2 - 2bm - 2\lambda & 2\bigl(-b + 2bm + \lambda\bigr)\\
    b - 2m             & b + 2m                       
  \end{bmatrix}
  \begin{bmatrix}
    \M\parentes{\Lambda - \frac{1}{2}, \frac{m - 1}{2}, \frac{b}{2}}\\
    \M\parentes{\Lambda + \frac{1}{2}, \frac{m - 1}{2}, \frac{b}{2}}
  \end{bmatrix}
  =
  \begin{bmatrix}
    0\\
    0
  \end{bmatrix}
\]
If it would hold that $\M\parentes{\Lambda - \frac{1}{2}, \frac{m - 1}{2},
\frac{b}{2}} = 0$, it would according to~\eqref{eq:eigenfunction} mean that
$\lambda$ is also an eigenvalue of the Dirichlet realization of $\Ham_{m-1}(b)$.
Since the smallest Dirichlet eigenvalue is greater than $b$, and we assume that
$\lambda < b$, this is not possible. Thus, the determinant of the coefficient
matrix must be zero,
\[
  \begin{aligned}
  0
  & =
  \det
  \begin{bmatrix}
    b^2 - 2bm - 2\lambda & 2\bigl(-b + 2bm+\lambda\bigr)\\
    b - 2m               & b + 2m
  \end{bmatrix}
  \\
  & =
  -4b\bracket*{\lambda - \parentes*{\frac{b}{2} - m}\parentes*{\frac{b}{2} - (m - 1)}}.
  \end{aligned}
\]
We conclude that if $\lambda < b$ is a common eigenvalue of $\Hm(b)$ and
$\Ham^{(m - 1)}(b)$, then
\[
  \lambda=\parentes*{\frac{b}{2} - m}\parentes*{\frac{b}{2} - (m - 1)}.\qedhere
\]
\end{proof}

The fact that it is the lowest eigenvalues of $\Ham^{(m - 1)}(b)$ and
$\Ham^{(m)}(b)$ that intersect follows from the upper bound $\lambda < b$. The
property in the theorem, however, seems also to be true for higher eigenvalues
of $\Ham^{(m - 1)}(b)$ and $\Ham^{(m)}(b)$. 

\subsection{The exterior of the disk}

For the exterior the situation is similar. We use again the same notation as
above, but with the operator now acting in \(L^2((1,+\infty),r\dd r)\). Here,
the other Whittaker function \(\W\) is used. Similar formulas exist for it. In
fact, the formula~\eqref{eq:diffM} is replaced by
\[
  \frac{\dd}{\dd z}\W(\kappa,\mu,z)
  =
  \left(\frac{1}{2}-\frac{\kappa}{z}\right)\W(\kappa,\mu,z)
  -\frac{1}{z}\W(\kappa+1,\mu,z),
\]
and then the formula~\eqref{eq:form} is replaced by
\[
  \left(\frac{\lambda}{b} + m + 1 - \frac{b}{2}\right)
  \W\left(\Lambda, \frac{m}{2}, \frac{b}{2}\right) 
  + 2\W\left(\Lambda + 1, \frac{m}{2}, \frac{b}{2}\right) = 0,
\]
and the formula~\eqref{eq:formm1} is replaced by
\[
  \left(\frac{\lambda}{b} + m - \frac{b}{2}\right)
  \W\left(\Lambda - \frac{1}{2}, \frac{m - 1}{2}, \frac{b}{2}\right)
  +2\W\left(\Lambda + \frac{1}{2}, \frac{m - 1}{2}, \frac{b}{2}\right) = 0.
\]
The formula~\eqref{eq:Mfinal} is replaced by
\[
  \W(\kappa, \mu, z) =   
  \frac{1}{\sqrt{z}}
  \left[
    \W\left(\kappa + \frac{1}{2}, \mu - \frac{1}{2}, z\right)
  + \left(\kappa + \mu - \frac{1}{2}\right) 
    \W\left(\kappa - \frac{1}{2}, \mu - \frac{1}{2}, z\right)
  \right],
\]
and the formula~\eqref{eq:Mfirstarg} is replaced by
\[
  \W(\kappa + 1, \mu, z) + (2\kappa - z) \W(\kappa, \mu, z) 
  + \left(\kappa - \mu - \frac{1}{2}\right)
    \left(\kappa + \mu - \frac{1}{2}\right) \W(\kappa - 1, \mu, z) = 0.
\]
Once that is done, the reasoning is similar, and in fact exactly the same
formula as for the interior of the disk holds.

\begin{proposition}\label{prop:intersections}
  Assume that $b > 0$ and $m \geq 2$. If $\lambda_1^{(m - 1)}(b) =
  \lambda_1^{(m)}(b) < b$ then $\lambda_1^{(m)}(b) = (b/2 - m)(b/2 - (m - 1))$.
\end{proposition}

\bibliographystyle{plain}
\bibliography{KLS.bib}

\end{document}